\documentclass[11pt]{article}
\usepackage{graphicx}
\usepackage{amssymb,amsfonts}
\usepackage{latexsym,bbm}
\usepackage{tikz,pgf}

\title{\textbf{The symmetric group action on rank-selected posets of injective words}}

\author{Christos~A.~Athanasiadis\\
Department of Mathematics\\
University of Athens\\
Athens 15784, Hellas (Greece)\\
\texttt{caath@math.uoa.gr}}

\date{\small October 28, 2016}

  \def\CC{{\mathbb C}}

  \def\ZZ{{\mathbb Z}}
  \def\RR{{\mathbb R}}
  \def\OO{{\mathbbm 1}}

  \def\tT{{\mathcal T}}
  \def\fS{\mathfrak{S}}

  \def\Des{{\rm Des}}

  \def\sm{\smallsetminus}

  \newcommand{\qed}{$\hfill \Box$}

\begin{document}
\maketitle

\newtheorem{theorem}{Theorem}[section]
\newtheorem{proposition}[theorem]{Proposition}
\newtheorem{corollary}[theorem]{Corollary}
\newtheorem{definition}[theorem]{Definition}
\newtheorem{remark}[theorem]{Remark}
\newtheorem{lemma}[theorem]{Lemma}
\newtheorem{example}[theorem]{Example}
\newtheorem{examples}[theorem]{Examples}
\newtheorem{conjecture}[theorem]{Conjecture}
\newtheorem{fact}[theorem]{Fact}
\newtheorem{question}[theorem]{Question}
\newtheorem{observation}[theorem]{Observation}
\newtheorem{claim}[theorem]{Claim}

\begin{abstract}
The symmetric group $\fS_n$ acts naturally on the poset
of injective words over the alphabet $\{1, 2,\dots,n\}$.
The induced representation on the homology of this poset
has been computed by Reiner and Webb. We generalize their
result by computing the representation of $\fS_n$ on the
homology of all rank-selected subposets, in the sense of
Stanley. A further generalization to the poset of
$r$-colored injective words is given.

\end{abstract}


\section{Introduction and results}
\label{sec:intro}

Let $P$ be a finite partially ordered set (poset, for
short) and $G$ be a subgroup of the group of
automorphisms of $P$. We assume that $P$ has a minimum
element $\hat{0}$ and a maximum element $\hat{1}$ and
that it is graded of rank $n+1$, with rank function
$\rho: P \to \{0, 1,\dots,n+1\}$ (basic definitions on
posets can be found in \cite[Chapter~3]{StaEC1}). Then
$G$ defines a permutation representation $\alpha_P(S)$,
induced by the action of $G$ on the set of maximal
chains of the rank-selected subposet
\begin{equation}
P_S \ = \ \{ x \in P: \, \rho(x) \in S \} \, \cup \,
\{\hat{0}, \hat{1}\}
\label{eq:P_S}
\end{equation}
for every $S \subseteq \{1, 2,\dots,n\}$ and one can
consider the virtual $G$-representation
\begin{equation}
\beta_P(S) \ = \ \sum_{T \subseteq S} \, (-1)^{|S-T|}
\, \alpha_P(T).
\label{eq:beta(S)}
\end{equation}
When $P$ is Cohen--Macaulay, $\beta_P(S)$ coincides
with the non-virtual $G$-representation induced on
the top homology group of (the order complex of)
$\bar{P}_S := P_S \sm \{\hat{0}, \hat{1}\}$.
The relations of the representations $\alpha_P(S)$
and $\beta_P(S)$ with one another and with the
structure of $P$ were investigated by Stanley in
his seminal work~\cite{Sta82}, where several
interesting examples (such as the symmetric group
action on Boolean, subspace and partition lattices
and the hyperoctahedral group action on the face
lattice of the cross-polytope) were analyzed.

This paper aims to analyze the symmetric group
action on another interesting poset, namely the
poset of injective words. We denote by $P_n$ the
set of words over the alphabet $\{1, 2,\dots,n\}$
with no repeated letter, partially ordered by
setting $u \le v$ if $u$ is a subword of $v$,
along with an artificial maximum element
$\hat{1}$ (the empty word is the minimum element
$\hat{0}$). The poset $P_n$ is the augmented
face poset of a Boolean regular cell complex;
see, for
instance,~\cite[Section~1]{ReiW04}. It was
introduced by Farmer~\cite{Fa79}, was shown to
be Cohen--Macaulay by Bj\"orner and
Wachs~\cite[Section~6]{BW83} and was further
studied, among other places, in~\cite{HH04,
JW09, KK13, ReiW04}.

\bigskip
\begin{center}
\begin{tikzpicture}[scale=1.4]
\label{fg:inj3}
   \draw(0,0) node(b){};
   \draw(-1,1) node(1){};
   \draw(0,1) node(2){};
   \draw(1,1) node(3){};
   \draw(-2.5,2) node(12){};
   \draw(-1.5,2) node(13){};
   \draw(-0.5,2) node(21){};
   \draw(0.5,2) node(23){};
   \draw(1.5,2) node(31){};
   \draw(2.5,2) node(32){};
   \draw(-2.5,3) node(123){};
   \draw(-1.5,3) node(132){};
   \draw(-0.5,3) node(213){};
   \draw(0.5,3) node(231){};
   \draw(1.5,3) node(312){};
   \draw(2.5,3) node(321){};
   \draw(0,4) node(t){};

   \draw(b) -- (1) -- (12) -- (123) -- (t);
   \draw(12) -- (132) -- (t);
   \draw(12) -- (312) -- (t);
   \draw(1) -- (21) -- (213) -- (t);
   \draw(21) -- (231) -- (t);
   \draw(21) -- (321) -- (t);
   \draw(1) -- (13) -- (123);
   \draw(13) -- (132);
   \draw(13) -- (213);
   \draw(1) -- (31) -- (231);
   \draw(31) -- (321);
   \draw(31) -- (312);
   \draw(b) -- (2) -- (23) -- (123);
   \draw(2) -- (12);
   \draw(2) -- (21);
   \draw(23) -- (213);
   \draw(23) -- (231);
   \draw(2) -- (32) -- (132);
   \draw(32) -- (312);
   \draw(32) -- (321);
   \draw(b) -- (3) -- (13) -- (123);
   \draw(13) -- (213);
   \draw(3) -- (32) -- (132);
   \draw(3) -- (31) -- (312);
   \draw(3) -- (23);
   \draw(31) -- (321);

   \node[circle,inner sep=0pt,fill=white] at (b) {$\hat{0}$};
   \node[circle,inner sep=0pt,fill=white] at (1){$1$};
   \node[circle,inner sep=0pt,fill=white] at (2){$2$};
   \node[circle,inner sep=0pt,fill=white] at (3){$3$};
   \node[circle,inner sep=0pt,fill=white] at (12){$12$};
   \node[circle,inner sep=0pt,fill=white] at (21){$21$};
   \node[circle,inner sep=0pt,fill=white] at (13){$13$};
   \node[circle,inner sep=0pt,fill=white] at (31){$31$};
   \node[circle,inner sep=0pt,fill=white] at (23){$23$};
   \node[circle,inner sep=0pt,fill=white] at (32){$32$};
   \node[circle,inner sep=0pt,fill=white] at (123){$123$};
   \node[circle,inner sep=0pt,fill=white] at (132){$132$};
   \node[circle,inner sep=0pt,fill=white] at (213){$213$};
   \node[circle,inner sep=0pt,fill=white] at (231){$231$};
   \node[circle,inner sep=0pt,fill=white] at (312){$312$};
   \node[circle,inner sep=0pt,fill=white] at (321){$321$};
   \node[circle,inner sep=0pt,fill=white] at (t){$\hat{1}$};

\end{tikzpicture}
\end{center}
\medskip

The symmetric group $\fS_n$ acts on $P_n$, by the
natural action of a permutation $w \in \fS_n$ on
each letter of a word over the alphabet $\{1,
2,\dots,n\}$, as a group of automorphisms. The
decomposition of the $\fS_n$-representation on the
top homology of $\bar{P}_n$ into irreducibles was
computed, using the Hopf trace formula, by
Reiner and Webb~\cite{ReiW04} and was further
refined by Hanlon and Hersh~\cite{HH04}, who also
discovered interesting connections with the spectrum
a of certain Markov chain on $\fS_n$; see \cite{DS15}
for recent developments. Recall that the set $\Des(Q)$
of descents of a standard Young tableau $Q$ of size $n$ consists of all integers $i \in \{1, 2,\dots,n-1\}$ for
which $i+1$ lies in a lower row that $i$. The result
of Reiner and Webb can be stated as follows.

\begin{theorem} {\rm (\cite[Proposition~2.3]{ReiW04})} \label{thm:ReiW04}
For $T = \{1, 2,\dots,n\}$, the multiplicity of the irreducible $\fS_n$-representation corresponding to
$\lambda \vdash n$ in $\beta_{P_n}(T)$ is equal to
the number of standard Young tableaux $Q$ of shape
$\lambda$ for which the smallest element of $\Des(Q)
\cup \{n\}$ is even.
\end{theorem}

To state our main result, we introduce the following
notation and terminology. Recall that the descent set
of a permutation $w \in \fS_n$ consists of
those integers $i \in \{1, 2,\dots,n-1\}$ for which
$w(i) > w(i+1)$. Let $S = \{s_1, s_2,\dots,s_k\}$
be a subset of $\{1, 2,\dots,n-1\}$ with $s_1 < s_2
< \cdots < s_k$ and denote by $w_S$ the element of
$\fS_n$ of largest possible Coxeter length (number
of inversions) with descent set equal to $S$. For
example, if $n = 8$ and $S = \{2, 3, 6\}$ then,
written in one-line notation, $w_S = (7, 8, 6, 3,
4, 5, 1, 2)$. For a permutation $w \in \fS_n$ with
descent set $S$ and standard Young tableaux $Q$ of
size $n$, let $\tau(w, Q)$ denote the largest
$i \in \{0, 1,\dots, k+1\}$ for which: (a) $w(x) =
w_S (x)$ for $x > s_{k-i+1}$; and (b) no descent of
$Q$ is smaller than $n - s_{k-i+1}$, where we have
set $s_0 = 0$ and $s_{k+1} = n$. Note that these
conditions are trivially satisfied for $i=0$ and that
for $S = \{1, 2,\dots,n-1\}$ (in which case $w = w_S$
is the only element of $\fS_n$ with descent set equal
to $S$), the integer $\tau(w, Q)$ is exactly the
smallest element of $\Des(Q) \cup \{n\}$. Thus, the
following statement generalizes
Theorem~\ref{thm:ReiW04}.

\begin{theorem} \label{thm:main}
For $T \subseteq \{1, 2,\dots,n\}$, denote by
$b_\lambda (T)$ the multiplicity of the irreducible
$\fS_n$-representation corresponding to $\lambda \vdash
n$ in $\beta_{P_n}(T)$. Then, for $S \subseteq \{1, 2,\dots,n-1\}$ and $\lambda \vdash n$,
\begin{itemize}
\itemsep=0pt
\item[{\rm (a)}]
$b_\lambda (S)$ is equal to the number of pairs $(w, Q)$
of permutations $w \in \fS_n$ with descent set $S$ and standard Young tableaux $Q$ of shape $\lambda$ for which $\tau(w, Q)$ is odd; and

\item[{\rm (b)}]
$b_\lambda (S \cup \{n\})$ is equal to the number
of pairs $(w, Q)$ of permutations $w \in \fS_n$
with descent set $S$ and standard Young tableaux
$Q$ of shape $\lambda$ for which $\tau(w, Q)$ is
even.
\end{itemize}
\end{theorem}

The structure and other results of this paper are
as follows. Section~\ref{sec:pre} fixes notation
and terminology. Section~\ref{sec:good} obtains a
formula (Equation~\ref{eq:lemgood}) and a partial
result (Equation~\ref{eq:bgood}) for the
representations $\beta_{P_n}(S)$ which are valid
for actions more general than that of the symmetric
group on the poset of injective words and discusses
examples. Section~\ref{sec:color} proves
Theorem~\ref{thm:main}, generalized to the action
of the symmetric group on the poset of $r$-colored
injective words, and discusses some consequences and
special cases.

\section{Preliminaries}
\label{sec:pre}

This section fixes some basic notation and
terminology which will be used in the sequel.
For all positive integers $n$, we set $[n] =:
\{1, 2,\dots,n\}$.

Our notation on posets follows mostly that of
\cite[Chapter~3]{StaEC1}, where basic background
on this topic can also be found. In particular,
we denote by $B_n$ the boolean lattice of all
subsets of the set $[n]$, partially ordered by
inclusion. A finite poset $P$ is \emph{graded}
of rank $n+1$
if every maximal chain in $P$ has length $n+1$
(meaning, cardinality $n+2$). The rank $\rho(x)$
of $x \in P$ is then defined as the length of the
longest chain in $P$ with maximum element equal
to $x$. Assume that such a poset $P$ has a minimum
and a maximum element. For $S \subseteq [n]$, we
denote by $a_P (S)$ the number of maximal chains
of the rank-selected subposet $P_S$ of $P$ defined
in (\ref{eq:P_S}). For instance, $a_P(\varnothing) =
1$ and $a_P([n])$ is the number of maximal chains
of $P$. The numbers $b_P(S)$ are defined by
setting
\begin{equation}
b_P(S) \ = \ \sum_{T \subseteq S} \, (-1)^{|S-T|}
\, a_P(T)
\label{eq:b(S)1}
\end{equation}
for $S \subseteq [n]$ or, equivalently, by
\begin{equation}
a_P(S) \ = \ \sum_{T \subseteq S} \, b_P(T)
\label{eq:b(S)2}
\end{equation}
for $S \subseteq [n]$.

An \emph{automorphism} of a poset $P$ is a bijective
map $\varphi: P \to P$ such that $x \le y
\Leftrightarrow \varphi(x) \le \varphi(y)$ for $x, y
\in P$. The automorphisms of $P$ form a group under
the law of composition of maps.

For basic background on the representation theory
of the symmetric group and the combinatorics of
Young tableaux, we refer the reader to~\cite{Sa01} \cite[Chapter~7]{StaEC2} \cite[Lecture~2]{Wa07}.
All representations we consider are defined over
the field of complex numbers $\CC$. The trivial
representation of a group $G$ is denoted by $\OO_G$.
We follow the English notation to describe Young
tableaux.

Let $P$ be a finite, graded poset of rank $n+1$
with a minimum and a maximum element and let $G$
be a subgroup of the group of automorphisms of
$P$. As discussed in the introduction, for $S
\subseteq [n]$, we denote by $\alpha_P(S)$ the
permutation $G$-representation induced by the
action of $G$ on the set of maximal chains of
$P_S$ and by $\beta_P(S)$ the virtual
$G$-representation defined by (\ref{eq:beta(S)})
or, equivalently, by
\begin{equation}
\alpha_P(S) \ = \ \sum_{T \subseteq S} \,
\beta_P(T)
\label{eq:beta(S)2}
\end{equation}
for $S \subseteq [n]$. We note that the dimension
of $\alpha_P(S)$ is equal to $a_P(S)$. If $P$ is
Cohen--Macaulay over $\CC$, then $\beta_P(S)$ is a
non-virtual $G$-representation
\cite[Theorem~1.2]{Sta82} whose dimension is equal
to $b_P(S)$; see \cite{Sta82} \cite[Section~3.4]{Wa07}
for more information on these representations.

\section{Good actions on simplicial posets}
\label{sec:good}

This section obtains a formula for the
representations $\beta_P(S)$ under certain
hypotheses, which are satisfied by the symmetric
group action on the poset of injective words.
Throughout this section, $P$ is a finite, graded
poset of rank $n+1$ with minimum element
$\hat{0}$ and maximum element $\hat{1}$ and $G$
is a subgroup of the group of automorphisms of $P$.

Following \cite[Section~3.16]{StaEC1}, we say that
$P \sm \{\hat{1}\}$ is \emph{simplicial} if the
closed interval $[\hat{0}, x]$ of $P$ is isomorphic
to a Boolean
lattice for every $x \in P \sm \{\hat{1}\}$. We
call the action of $G$ on $P$ \emph{good} if $g
\cdot y = y$ for some $g \in G$ and $y \in P \sm
\{\hat{1}\}$ implies that $g \cdot x = x$ for
every $x \in P$ with $x \le y$. Let us also denote
by $b_n(S)$ the number of permutations $w \in \fS_n$
with descent set $S$.
\begin{theorem} \label{thm:good}
Suppose that the poset $P \sm \{\hat{1}\}$ is
simplicial and that the action of $G$ is good. Then,
\begin{equation} \label{eq:lemgood}
\beta_P (S) \ = \ (-1)^k \, \OO_G \ + \ \sum_{i=1}^k
\, (-1)^{k-i} \, b_{s_i} (\{s_1, s_2,\dots,s_{i-1}\})
\cdot \alpha_P (\{s_i\})
\end{equation}
for every $S = \{s_1, s_2,\dots,s_k\} \subseteq [n]$
with $s_1 < s_2 < \cdots < s_k$. In particular,
\begin{equation} \label{eq:bgood}
\beta_P (S) \, + \, \beta_P (S \cup \{n\}) \ = \
b_n(S) \cdot \alpha_P(\{n\})
\end{equation}
for every $S \subseteq [n-1]$.
\end{theorem}

\noindent
\emph{Proof}. We set $a_n (T) := a_{B_n} (T)$ for
$T \subseteq [n-1]$. Thus, $a_n (T)$ is the
multinomial coefficient equal to the number of
chains $c$ in the Boolean lattice $B_n$ for which
the set of cardinalities (ranks) of the elements
of $c$ is equal to $T$. It is well known
\cite[Corollary~3.13.2]{StaEC1} that
\begin{equation}
\sum_{T \subseteq S} \, (-1)^{|S-T|} \, a_n(T) \
= \ b_n(S)
\label{eq:bn(S)}
\end{equation}
for every $S \subseteq [n-1]$, where $b_n(S)$ is as
in the sentence preceding the theorem.

We claim that
\begin{equation}
\alpha_P(S) \ = \ a_m (S \sm \{m\}) \cdot
\alpha_P(\{m\})
\label{eq:agood}
\end{equation}
for every $S \subseteq [n]$ with maximum element
equal to $m$. Indeed, since the action of $G$ on
the poset $P$ is good, a maximal chain $c$ of $P_S$
is fixed by an element $g \in G$ if and only if
the maximum element of $c \sm \{\hat{1}\}$ is
fixed by $g$. As a result, and since $P \sm
\{\hat{1}\}$ is simplicial, the number of maximal
chains of $P_S$ fixed by $g$ is equal to the
product of the number of
elements of $P$ of rank $m$ which are fixed by
$g$ with the number $a_m (S \sm \{m\})$ of chains
in the Boolean lattice $B_m$ whose
elements have set of ranks equal to $S$. This
shows that the characters of the representations
in the two hand sides of (\ref{eq:agood}) are equal
and verifies our claim.

Using (\ref{eq:agood}), for $S = \{s_1,
s_2,\dots,s_k\} \subseteq [n]$ with $s_1 < s_2 <
\cdots < s_k$, the defining equation~(\ref{eq:beta(S)})
gives
\begin{eqnarray*}
\beta_P(S) & = & (-1)^k \, \OO_G \ \, + \ \
\sum_{i=1}^k \! \! \!
\sum_{\scriptsize \begin{array}{c} T \subseteq S:
\\ \max(T) = s_i \end{array}} \! \! (-1)^{|S-T|} \,
a_{s_i} (T \sm \{s_i\}) \cdot \alpha_P (\{s_i\}) \\
& & \\
& = & (-1)^k \, \OO_G \ \, + \ \ \sum_{i=1}^k \
\sum_{T \subseteq \{s_1,\dots,s_{i-1}\}} \
(-1)^{k-|T|-1} \, a_{s_i} (T) \cdot \alpha_P
(\{s_i\}).
\end{eqnarray*}
Using~(\ref{eq:bn(S)}) to compute the inner sum
gives the expression (\ref{eq:lemgood}) for
$\beta_P(S)$. This expression directly implies
Equation~(\ref{eq:bgood}). Alternatively, using
the defining equation (\ref{eq:beta(S)}) we find
that, for $S \subseteq [n-1]$,
\begin{eqnarray*}
\beta_P (S \cup \{n\}) & = & \sum_{T \subseteq S}
\, (-1)^{|S-T|} \, \alpha_P (T \cup \{n\}) \ - \
\sum_{T \subseteq S} \, (-1)^{|S-T|} \, \alpha_P (T)
\\ & = & \sum_{T \subseteq S} \, (-1)^{|S-T|} \,
\alpha_P (T \cup \{n\}) \ - \ \beta_P (S).
\end{eqnarray*}
Using (\ref{eq:bn(S)}) and (\ref{eq:agood}), we
conclude that
\begin{eqnarray*}
\beta_P (S) \, + \, \beta_P (S \cup \{n\}) & = &
\sum_{T \subseteq S} \, (-1)^{|S-T|} \, \alpha_P
(T \cup \{n\}) \ = \ \sum_{T \subseteq S} \,
(-1)^{|S-T|} \, a_n (T) \cdot \alpha_P (\{n\}) \\
& = & b_n (S) \cdot \alpha_P (\{n\})
\end{eqnarray*}
and the proof follows.
\qed

\begin{example} \rm
Let $n, r$ be positive integers and let $\ZZ_r$
denote the abelian group of integers modulo $r$.
We consider the
cartesian product $[n] \times \ZZ_r$ as an alphabet
and think of the second coordinate of an element
of this set as one of $r$ possible colors attached
to the first coordinate. We denote by $P_{n,r}$ the
set of words over this alphabet whose letters have
pairwise distinct first coordinates, partially
ordered by the subword order (so that $u \le v$
if and only if $u$ can be obtained by deleting
some of the letters of $v$), along with an
artificial maximum element $\hat{1}$ (the empty
word is the minimum element $\hat{0}$). This 
is the \emph{poset of $r$-colored injective
words}; it specializes to the poset of injective
words $P_n$, discussed in the introduction, when
$r=1$. The Cohen-Macaulayness of $P_{n,r}$ follows
from \cite[Theorem~1.2]{JW09}.

The symmetric group $\fS_n$ acts on $P = P_{n,r}$
by acting on the first coordinate of each letter of
a word and leaving the colors unchanged. We leave
it to the reader to verify that the hypotheses of
Theorem~\ref{thm:good} are satisfied for this
action. Clearly, $\alpha_P(\{n\})$ is isomorphic
to the direct sum of $r^n$ copies of the regular representation $\rho^{\rm reg}$ of $\fS_n$ and
thus, Equation~(\ref{eq:bgood}) gives
\begin{equation} \label{eq:betacolor}
\beta_P (S) \, + \, \beta_P (S \cup \{n\}) \ = \
b_n(S) \, r^n \cdot \rho^{\rm reg}
\end{equation}
for every $S \subseteq [n-1]$. Equivalently,
denoting by $b_{r, \lambda} (T)$ the multiplicity
of the irreducible $\fS_n$-representation
corresponding to $\lambda \vdash n$ in $\beta_P
(T)$, we have
\begin{equation} \label{eq:blambdacolor}
b_{r, \lambda} (S) \, + \, b_{r, \lambda} (S \cup
\{n\}) \ = \ b_n(S) \, r^n f^\lambda
\end{equation}
for every $\lambda \vdash n$, where $f^\lambda$
is the number of standard Young tableaux of shape
$\lambda$. Equation~(\ref{eq:blambdacolor}) suggests
that there is a combinatorial interpretation to
each of the summands on its left-hand side in terms
of pairs of elements of $\fS_n$ with descent set $S$
and $r$-colored standard Young tableaux of shape
$\lambda$. Theorem~\ref{thm:cmain}, proved in the
following section, provides such an interpretation
and thus determines the corresponding summands
of~(\ref{eq:betacolor}).

We note that Equation~(\ref{eq:betacolor})
substantially refines
\begin{equation} \label{eq:betacolortrivial}
\sum_{T \subseteq [n]} \beta_P (T) \ = \
r^n \, n! \cdot \rho^{\rm reg},
\end{equation}
which can be shown directly by observing that the
permutation representation $\alpha_P([n])$ of
$\fS_n$ on the set of maximal chains of $P$ is
isomorphic to the direct sum of $r^n \, n!$ copies
of the regular representation of $\fS_n$.
\end{example}

\begin{remark} \rm
There are larger groups than $\fS_n$, such as the 
wreath products $\fS_n[\ZZ_r]$ and $\fS_n[\fS_r]$,
which act naturally on the poset $P_{n,r}$ of 
$r$-colored injective words by good actions. We 
leave it as an open problem to refine 
Theorem~\ref{thm:cmain} in this direction and note 
that the induced $\fS_n[\ZZ_2]$-representation on 
the top homology group of $\bar{P}_{n,2}$ was 
computed in~\cite[Theorem~7.3]{AAER15}.  
\qed
\end{remark}

\begin{example} \rm
Let $\Delta$ be a finite, abstract simplicial complex
of dimension $n-1$. We assume that $\Delta$ is pure
(every maximal face of $\Delta$ has dimension $n-1$)
and (completely) balanced (the vertices of $\Delta$
are colored with $n$ colors, so that vertices of any
maximal face of $\Delta$ have distinct colors; see
\cite[Section~III.4]{StaCCA} for more information
about this class of complexes).

Let $P_\Delta$ be the set of faces of $\Delta$,
ordered by inclusion, with an artificial maximum
element $\hat{1}$ attached, and let $G$ be a subgroup
of the group of automorphisms of $P_\Delta$ whose
action preserves the colors of the vertices of
$\Delta$. We leave it again to the reader to verify
that the hypotheses of Theorem~\ref{thm:good}
are satisfied for the action of $G$ on $P_\Delta$.
The special case in which $\Delta$ is the order
complex of a poset was analysed by
Stanley~\cite[Section~8]{Sta82}.
\end{example}

\section{The poset of colored injective words}
\label{sec:color}

This section applies Theorem~\ref{thm:good} to
decompose the $\fS_n$-representations $\beta_P (S)$
into irreducibles for the poset of $r$-colored
injective words, thus generalizing
Theorem~\ref{thm:main}.

To state the main result of this section, we need
to modify some of the definitions introduced before
Theorem~\ref{thm:main}. An \emph{$r$-colored standard
Young tableau} of shape $\lambda \vdash n$ is a
standard Young tableau of shape $\lambda$, each
entry of which has been colored with one of the
elements of $\ZZ_r$. Let $S = \{s_1, s_2,\dots,s_k\}
\subseteq [n-1]$ with $s_1 < s_2 < \cdots < s_k$ and
recall the definition of $w_S \in \fS_n$, given in
the introduction. Given a permutation $w \in \fS_n$
with descent set $S$ and an $r$-colored standard
Young tableau $Q$ of size $n$, we denote by
$\tau(w, Q)$ the largest $i \in \{0, 1,\dots,k+1\}$
for which: (a) $w(x) = w_S (x)$ for $x > s_{k-i+1}$;
and (b) all numbers $1, 2,\dots,n - s_{k-i+1}$ appear
in the first row of $Q$ and are colored with the zero
color, where we have set $s_0 = 0$ and $s_{k+1} = n$.
Note that these conditions are vacuously satisfied
for $i=0$ and that the definition of $\tau(w, Q)$
agrees with that given in the introduction in the
special case $r=1$.

\begin{theorem} \label{thm:cmain}
Let $P$ be the poset $P_{n,r}$ of $r$-colored
injective words and for $T \subseteq [n]$, denote
by $b_{r, \lambda} (T)$ the multiplicity of the
irreducible $\fS_n$-representation corresponding
to $\lambda \vdash n$ in $\beta_P(T)$. Then, for
$S \subseteq [n-1]$ and $\lambda \vdash n$,
\begin{itemize}
\itemsep=0pt
\item[{\rm (a)}]
$b_{r, \lambda} (S)$ is equal to the number of pairs
$(w, Q)$ of permutations $w \in \fS_n$ with descent
set $S$ and $r$-colored standard Young tableaux $Q$
of shape $\lambda$ for which $\tau(w, Q)$ is odd;
and

\item[{\rm (b)}]
$b_{r, \lambda} (S \cup \{n\})$ is equal to the number
of pairs $(w, Q)$ of permutations $w \in \fS_n$
with descent set $S$ and $r$-colored standard Young
tableaux $Q$ of shape $\lambda$ for which $\tau(w,
Q)$ is even.
\end{itemize}
\end{theorem}

\noindent
\emph{Proof}.
To apply Theorem~\ref{thm:good} and~(\ref{eq:agood}),
we need to determine $\alpha_P (\{m\})$ for $m \in
[n]$. The elements of $P$ of rank $m$ are the
$r$-colored injective words of length $m$ over the
alphabet $[n]$. Clearly, the action of $\fS_n$ on
these words has $r^m$ orbits, corresponding to the
$r^m$ possible coloring patterns, and the stabilizer
of the action on each orbit is isomorphic to the
Young subgroup $(\fS_1)^m \times \fS_{n-m}$ of
$\fS_n$. Using Young's
rule~\cite[Theorem~2.12.2]{Sa01} to decompose the
permutation representation of $\fS_n$ on the set of
left cosets of this subgroup we find that
\begin{equation} \label{eq:s}
\alpha_P (\{m\}) \ = \ r^m \ \sum_{\lambda \vdash
n} \, f^{\lambda, n-m} \rho^\lambda,
\end{equation}
where $\rho^\lambda$ stands for the irreducible
$\fS_n$-representation corresponding to $\lambda
\vdash n$ and $f^{\lambda, n-m}$ denotes the number
of standard Young tableaux of size $n$ whose first
row contains the numbers $1, 2,\dots,n-m$. Now let
$S = \{s_1, s_2,\dots,s_k\} \subseteq [n-1]$ with
$s_1 < s_2 < \cdots < s_k$ and set $s_0 := 0$ and
$s_{k+1} := n$. Equation~(\ref{eq:lemgood}) implies
the expression
\begin{equation} \label{eq:blS1}
b_{r, \lambda} (S) \ = \ (-1)^k \, \delta_{\lambda,
(n)} \ + \ \sum_{i=1}^k \, (-1)^{k-i} \, b_{s_i}
(\{s_1, s_2,\dots,s_{i-1}\}) \, r^{s_i}
f^{\lambda, n-s_i}
\end{equation}
for the multiplicity of $\rho^\lambda$ in $\beta_P(S)$,
where $\delta_{\lambda, (n)}$ is a Kronecker delta.
For $1 \le i \le k$, clearly, $r^{s_i} f^{\lambda,
n-s_i}$ is equal to the number of $r$-colored standard
Young tableaux of shape $\lambda$ whose first row
contains the numbers $1, 2,\dots,n-s_i$, all colored
with the zero color. We may also interpret $b_{s_i}
(\{s_1, s_2,\dots,s_{i-1}\})$ as the number of
permutations $w \in \fS_n$ with descent set $S$ for
which $w(x) = w_S (x)$ for $x > s_i$. Therefore, Equation~(\ref{eq:blS1}) can be rewritten as
\begin{equation} \label{eq:blS2}
b_{r, \lambda} (S) \ = \ \sum_{i=1}^{k+1} \,
(-1)^{i-1} \, |\tT_i|,
\end{equation}
where $\tT_i$ is the set of all pairs $(w, Q)$ of
permutations $w \in \fS_n$ with descent set $S$ and
$r$-colored standard Young tableaux $Q$ of shape
$\lambda$ for which: (a) $w(x) = w_S (x)$ for $x >
s_{k-i+1}$; and (b) all numbers $1, 2,\dots,n -
s_{k-i+1}$ appear in the first row of $Q$ and are
colored with the zero color. Noting that $\tT_0$ is
the set of all pairs $(w, Q)$ of permutations $w
\in \fS_n$ with descent set $S$ and $r$-colored
standard Young tableaux $Q$ of shape $\lambda$ and
that $\tT_0 \supseteq \tT_1 \supseteq \tT_2
\supseteq \cdots$ we conclude that
\smallskip
\begin{eqnarray}
b_{r, \lambda} (S) &=& |\tT_1 \sm \tT_2| \, +
\, |\tT_3 \sm \tT_4| \, + \, \cdots
\label{eq:blS3} \\
b_{r, \lambda} (S \cup \{n\}) &=& |\tT_0 \sm \tT_1|
\, + \, |\tT_2 \sm \tT_3| \, + \, \cdots,
\label{eq:blS4}
\end{eqnarray}

\noindent
where the first equation follows from (\ref{eq:blS2})
and the second by a similar argument, or from the
first by appealing to the last statement of
Theorem~\ref{thm:good}. The two equations are
equivalent to parts (a) and (b), respectively, of
the theorem.
\qed

\begin{example} \rm
Keeping the notation of Theorem~\ref{thm:cmain},
let us consider the multiplicity of the trivial and
the sign representation in $\beta_P(T)$ for
$T\subseteq [n]$, corresponding to $\lambda = (n)$
and $\lambda = (1^n)$, respectively.
Equation~(\ref{eq:blambdacolor}) holds with
$f^\lambda = 1$ in both cases, for $S \subseteq
[n-1]$. For $\lambda = (1^n)$, the tableaux $Q$ in
the definition of $\tT_i$ have a single column and
for $\tT_i$ to be nonempty, we must have either
$i=0$, or $i = 1$ and $n-1 \in S$. From
(\ref{eq:blS3}) and (\ref{eq:blS4}) we infer that
\[ b_{r, (1^n)} (S) \ = \ \cases{
        0, & if \ $n-1 \not\in S$, \cr
r^{n-1} \, b_{n-1}(S \sm \{n-1\}), & if \ $n-1 \in
S$} \]
and
\[ b_{r, (1^n)} (S \cup \{n\}) \ = \ \cases{
     r^n \, b_n(S), & if \ $n-1 \not\in S$, \cr
     r^n \, b_n(S) - r^{n-1} \, b_{n-1}(S \sm \{n-1\}),
     & if \ $n-1 \in S$.} \]
Similarly, for $\lambda = (n)$, the pairs $(w, Q)$ of
permutations $w \in \fS_n$ and $r$-colored standard
Young tableaux $Q$ of shape $\lambda$ may be
identified with the elements of the wreath product 
$\fS_n[\ZZ_r]$, viewed as the $r$-colored permutations 
of the set $[n]$. Thus, Theorem~\ref{thm:cmain} gives 
a combinatorial interpretation to the multiplicity 
$b_{r, (n)} (T)$ in terms of such permutations and 
their descent and color patterns.
\end{example}

\begin{example} \rm
Suppose $S = [n-1]$. Then, $w_S$ is the only element
of $\fS_n$ with descent set equal to $S$. Therefore,
for $\lambda \vdash n$, the multiplicity $b_{r, \lambda} ([n-1])$ (respectively, $b_{r, \lambda} ([n])$ is
equal to the number of $r$-colored standard Young
tableaux $Q$ of shape $\lambda$ for which the largest
$i \in \{0, 1,\dots,n\}$ such that all numbers $1,
2,\dots,i$ appear in the first row of $Q$ and are
colored with the zero color is odd (respectively,
even). This extends Theorem~\ref{thm:ReiW04} to
$r$-colored injective words.
\end{example}

\begin{remark} \rm
As discussed in the beginning of
the proof of Theorem~\ref{thm:cmain}, we have $$
\alpha_P (\{m\}) \ = \ r^m \, \cdot \, \OO \uparrow^{\, \fS_n}_{\, (\fS_1)^m \times \fS_{n-m}} $$ for $n \ge 
m$, where $P = P_{n,r}$. By~\cite[Theorem~2.8]{Ch12}
\cite[Lemma~2.2]{RH15}, this expression and
(\ref{eq:agood}) show that for a fixed set $S$ of
positive integers, the representation $\beta_P(S)$
stabilizes for $n \ge 2 \max(S)$, in the sense
of~\cite{CF13}.
\qed
\end{remark}

As before, we view the elements of the group $\fS_n
[\ZZ_r]$ as the $r$-colored permutations of the set 
$[n]$, meaning permutations $u$ of $[n]$ with the 
entries $u(i)$ in their one-line notation colored, 
each with one of the elements of $\ZZ_r$. Let $S = 
\{s_1, s_2,\dots,s_k\} \subseteq [n-1]$ be as in the 
paragraph preceding Theorem~\ref{thm:cmain}. 
Given a permutation $w \in \fS_n$ with descent set 
$S$ and an $r$-colored permutation $u \in \fS_n
[\ZZ_r]$, we denote by $\tau(w, u)$ the largest $i 
\in \{0, 1,\dots,k+1\}$ for which: 
(a) $w(x) = w_S (x)$ for $x > s_{k-i+1}$; and (b) the 
numbers $u(1), u(2),\dots,u(n - s_{k-i+1})$ are all
colored with the zero color in $u$ and are in
increasing order.

\begin{corollary} \label{cor:cmain}
Let $P$ be as in Theorem~\ref{thm:cmain}. Then,
for $S \subseteq [n-1]$,
\begin{itemize}
\itemsep=0pt
\item[{\rm (a)}]
$b_P (S)$ is equal to the number of pairs $(w, u)$
of permutations $w \in \fS_n$ with descent set $S$
and $r$-colored permutations $u \in \fS_n[\ZZ_r]$ 
for which $\tau(w, u)$ is odd; and

\item[{\rm (b)}]
$b_P (S \cup \{n\})$ is equal to the number of pairs
$(w, u)$ of permutations $w \in \fS_n$ with descent
set $S$ and $r$-colored permutations $u \in \fS_n
[\ZZ_r]$ for which $\tau(w, u)$ is even.
\end{itemize}
\end{corollary}

\noindent
\emph{Proof}.
Since $b_P(T)$ is the dimension of the
$\fS_n$-representation $\beta_P(T) = \sum_{\lambda
\vdash n} b_{r, \lambda} (T) \rho^\lambda$, we have
\begin{equation} \label{eq:br(T)}
b_P(T) \ = \ \sum_{\lambda \vdash n} b_{r, \lambda} (T)
f^\lambda
\end{equation}
for $T \subseteq [n]$. This equation, combined with
part (a) of Theorem~\ref{thm:cmain}, implies that $b_P
(S)$ counts the triples $(w, P_0, Q_0)$ of permutations
$w \in \fS_n$ with descent set $S$, standard Young
tableaux $P_0$ of size $n$ and $r$-colored standard
Young tableaux $Q_0$ of the same shape as $P_0$ for
which $\tau(w, Q_0)$ is odd. A trivial extension of
the Robinson-Schensted correspondence (see, for
instance,~\cite[Section~3.1]{Sa01}) provides a bijection
from the group $\fS_n[\ZZ_r]$ of $r$-colored permutations 
to the set of pairs $(P_0, Q_0)$ of standard Young tableaux
$P_0$ of size $n$ and $r$-colored standard Young tableaux
$Q_0$ of the same shape as $P_0$ (under this
correspondence, the entry $i \in [n]$ of the tableau
$Q_0$ associated to $u \in \fS_n[\ZZ_r]$ is colored with
the color that $u(i)$ has in $u$). Standard properties
of the Robinson-Schensted correspondence imply that
$\tau(w,u) = \tau(w, Q_0)$, if $u \in \fS_n[\ZZ_r]$ is
mapped to $(P_0, Q_0)$ under the extended
correspondence. Thus, the set of triples $(w, P_0,
Q_0)$ mentioned above, the cardinality of which is
equal to $b_P(S)$, bijects to the set of pairs
$(w, u)$ mentioned in part (a). A similar argument
works for part (b).
\qed

\begin{remark} \rm
Part (b) of Corollary~\ref{cor:cmain} implies that
$b_P([n])$ is equal to the number, say $E_{n,r}$,
of $r$-colored permutations $u \in \fS_n[\ZZ_r]$ 
for which the largest $i \in \{0, 1,\dots,n\}$ such
that the numbers $u(1), u(2),\dots,u(i)$ are in
increasing order and all colored with the zero
color is even. On the other hand, from
Equation~(\ref{eq:lemgood}) we get
\begin{eqnarray*}
b_P([n]) & = & \dim \beta_P([n]) \ = \ (-1)^n \ +
\ \sum_{i=1}^n \, (-1)^{n-i} \, a_P(\{i\}) \\
& = & \sum_{i=1}^n \, (-1)^{n-i} \, r^i \frac{n!}
{(n-i)!} \ = \ D_{n,r},
\end{eqnarray*}
where $D_{n,r}$ is the number of derangements
(elements without fixed points of zero color) in
$\fS_n[\ZZ_r]$. The fact that $D_{n,r} = E_{n,r}$ was 
first discovered for $r=1$ by D\'esarm\'enien~\cite{De84};
for a much stronger statement, see
\cite[Theorem~1]{DeWa88} \cite[Corollary~3.3]{DeWa93}
for $r=1$ and \cite[Theorem~7.3]{AAER15} for $r=2$. In
fact, this stronger statement can be expressed in terms
of the Frobenius characteristic of the
$\fS_n[\ZZ_r]$-representation on the homology of
$\bar{P}_{n,r}$, induced by the action of $\fS_n
[\ZZ_r]$ on $P_{n,r}$; see \cite[Section~2]{ReiW04}
for $r=1$ and \cite[Theorem~7.3]{AAER15} for $r=2$.
\end{remark}

\end{document}